\documentclass[a4paper, BCOR=0.0mm, DIV=calc]{scrartcl}
\usepackage[T1]{fontenc}
\usepackage[utf8]{inputenc}
\usepackage[ngerman]{babel}
\usepackage[osf,sc]{mathpazo}
\usepackage{textcomp}
\usepackage{microtype}
\usepackage{amsmath, amssymb, amsfonts}
\usepackage{tabularx}
\usepackage{nicefrac}
\KOMAoptions{twoside=false, twocolumn=false, headinclude=false, footinclude=false, mpinclude=false, pagesize=auto}
\recalctypearea

\setkomafont{section}{\mdseries\scshape\Large}
\setkomafont{title}{\mdseries\scshape}

\begin{document}
\title{Eine kurze Methode, Summen unendlicher Reihen durch Differentialformeln zu untersuchen\footnote{
Originaltitel: "`Methodus succincta summas serierum infinitarum per formulas differentiales investigandi"', erstmals publiziert in "`\textit{Memoires de l'academie des sciences de St.-Petersbourg} 5, 1815, pp. 45-56"', Nachdruck in "`\textit{Opera Omnia}: Series 1, Volume 16, pp. 200 - 213"', Eneström-Nummer E746, übersetzt von: Alexander Aycock, Textsatz: Artur Diener,  im Rahmen des Projektes "`Eulerkreis Mainz"' }}
\author{Leonhard Euler}
\date{}
\maketitle
\paragraph{§1}
Auch wenn ich diesen Gegenstand schon des Öfteren betrachtet habe, sind die meisten Dinge, die sich darauf beziehen, Summen angenehm auszudrücken, durch verschiedene Bücher hindurch verteilt und auch durch Umwege gefunden; deshalb werde ich hier eine kurze Methode angeben, mit deren Hilfe die Summe einer bestimmten Reihe mit leichter Rechnung ohne Umwege durch eine sehr einfache Form untersucht werden können wird.
\paragraph{§2}
Es sei also $X$ eine bestimmte Funktion von $x$, und $X'$, $X''$, $X'''$, etc mögen daher entstehen, wenn anstelle von $x$ nacheinander $x+1$, $x+2$, $x+3$, etc geschrieben wird. Daher werden also jene Buchstaben $X$, $X'$, $X''$, $X'''$, etc für mich die Terme, die den Indizes einer gewissen Reihe $x$, $x+1$, $x+2$, $x+3$, etc entsprechen, bezeichnen. Nachdem diese Dinge festgesetzt worden sind, werde ich zwei Fälle unendlicher Reihen betrachten, von welchen im ersten alle Terme mit demselben Zeichen "`$+$"' versehen fortschreiten, sodass die zu summierende Reihe
\[
	X + X' + X'' + X''' + \mathrm{etc}
\]
ist. In dem anderen Fall aber mögen die Terme mit alternierenden Zeichen fortschreiten, sodass die zu summierende Reihe
\[
	X - X' + X'' - X''' + \mathrm{etc} 
\]
ist. Ich werde also diese zwei Fälle getrennt voneinander entwickeln.
\section*{Fall 1 \\ Summation der unendlichen Reihe \\ $S = X + X' + X'' + X''' + \mathrm{etc.}$}
Es bezeichne $S'$ die Summe derselben Reihe, welcher der erste Term entfernt wurde, sodass
\[
	S' = X' + X'' + X''' + \mathrm{etc}
\]
und es wird, weil $S$ eine gewisse Funktion von $x$ ist, die wir hier hauptsächlich untersuchen, $S'$ die gleiche Funktion von $x+1$ sein. Es ist also klar, dass $S-S' = X$ sein wird. Weil daher
\[
	S' = S + \partial S + \tfrac{1}{2}\partial\partial S + \mathrm{etc}
\]
ist, wo ich die Nenner, die die Potenzen des Elements $\partial x$ enthalten, um für Kürze zu sorgen, weil man sie sich ja von selbst merken kann, weglasse, so wird daher unsere Gleichung diese Form annehmen:
\[
	0 = X + \partial S + \tfrac{1}{2}\partial\partial S + \tfrac{1}{6}\partial^3 S + \tfrac{1}{24}\partial^4 S + \mathrm{etc.}
\]
\paragraph{§4}
Wenn daher also diese Reihe sehr stark konvergiert, wird beinahe $\partial S = -X$ sein und daher $S = -\int X \partial x$, welches Integral durch eine Konstante so zu bestimmen ist, dass es für ein unendlich groß genommenes $x$ verschwindet, deshalb weil die infinitesimalen Terme für $0$ gehalten werden können, weil andernfalls die Reihe selbst keine endliche Summe hätte. Nachdem die Summe näherungsweise bekannt ist, wollen wir für die wahre Summe
\[
	S = - \int X\partial x - \alpha X - \beta \partial X - \gamma \partial\partial X - \mathrm{etc}
\]
setzen und es wird daher
\[
	\partial S = -X - \alpha \partial X - \beta \partial\partial X - \gamma \partial^3 X - \mathrm{etc}
\]
sein. Wenn daher gleich für die einzelnen Differentiale von $S$ die daher entstehenden Werte eingesetzt werden, wird man zur folgenden Gleichung gelangen:
\[
\left.
\begin{array}{ccccccccccccc}
&+& X &-& \alpha \partial X &-& \beta \partial\partial X &-& \gamma\partial^3 X &-& \delta\partial^4 X &-& \mathrm{etc.} \\[1mm]
&-& X &-& \tfrac{1}{2} &-& \tfrac{1}{2}\alpha &-& \tfrac{1}{2}\beta &-& \tfrac{1}{2}\gamma &-& \mathrm{etc.} \\[1mm]
& & & & &-& \tfrac{1}{6} &-& \tfrac{1}{6}\alpha &-& \tfrac{1}{6}\beta &-& \mathrm{etc.} \\[1mm]
& & & & & & &-& \tfrac{1}{24} &-& \tfrac{1}{24}\alpha &-& \mathrm{etc.} \\[1mm]
& & & & & & & & &-& \tfrac{1}{120} &-& \mathrm{etc.}                  
\end{array}
\right\} = 0
\]
und die unbekannten Koeffizienten $\alpha$, $\beta$, $\gamma$, etc müssen gleich aus den folgenden Gleichungen bestimmt werden:
\[
	\alpha + \tfrac{1}{2} = 0,\quad \beta + \tfrac{1}{2}\alpha + \tfrac{1}{6} = 0,\quad \gamma + \tfrac{1}{2}\beta + \tfrac{1}{6}\alpha + \tfrac{1}{24} = 0,
\]
woher
\[
	\alpha = -\tfrac{1}{2},\quad \beta = \tfrac{1}{12},\quad \gamma = 0,\quad \mathrm{etc}
\]
wird.
\paragraph{§5}
Auf diese Weise aber würde das Finden der Buchstaben $\alpha$, $\beta$, $\gamma$, etc zu aufwändig werden, und trotzdem würde man nicht ein einziges Bildungsgesetz erkennen, nach welchen sie weiter fortschreiten könnten; deshalb werde ich auf eine völlig einzigartige Weise die Werte dieser Buchstaben untersuchen. Ich werde hier natürlich eine gewöhnliche Reihe betrachten, die nach denselben Koeffizienten fortschreitet, welche
\[
	V = 1 + \alpha z + \beta z^2 + \gamma z^3 + \delta z^4 + \mathrm{etc}
\]
sei, und es ist klar, wenn die Summe $V$ dieser Reihe auf eine endliche Form gebracht werden kann, dass dann, wenn dieselbe nach Potenzen von $z$ entwickelt wird, notwendigerweise dieselbe Reihe hervorgehen muss, nach welcher Übereinkunft die Werte der Buchstaben $\alpha$, $\beta$, $\gamma$, $\delta$, etc von selbst bekannt werden.
\paragraph{§6}
Aus diesen Relationen, die zwischen den Buchstaben $\alpha$, $\beta$, $\gamma$, $\delta$ etc bestehen und oben in §$4$ erwähnt worden sind, führe man die folgenden Operationen durch:
\[
\begin{array}{rccccccccccccccccc}
	V &=& 1 &+& \alpha z &+& \beta z^2 &+& \gamma z^3 &+& \delta z^4 &+& \varepsilon z^5 &+& \mathrm{etc.} \\[1mm]
	\tfrac{1}{2}zV &=& &+& \tfrac{1}{2} &+& \tfrac{1}{2}\alpha &+& \tfrac{1}{2}\beta &+& \tfrac{1}{2}\gamma &+& \tfrac{1}{2}\delta &+& \mathrm{etc.}\\[1mm]
	\tfrac{1}{6}zzV &=& & & &+& \tfrac{1}{6} &+& \tfrac{1}{6}\alpha &+& \tfrac{1}{6}\beta &+& \tfrac{1}{6}\gamma &+& \mathrm{etc.} \\[1mm]
	\tfrac{1}{24}z^3V &=& & & & & &+& \tfrac{1}{24} &+& \tfrac{1}{24}\alpha &+& \tfrac{1}{24}\beta &+& \mathrm{etc.} \\[1mm]
	\tfrac{1}{120}z^4V &=& & & & & & & &+& \tfrac{1}{120} &+& \tfrac{1}{120}\alpha &+& \mathrm{etc.} \\[1mm]
	\tfrac{1}{720}z^5V &=& & & & & & & & & &+& \tfrac{1}{720} &+& \mathrm{etc.} \\[1mm]
	\mathrm{etc.}
\end{array}
\]
Auf diese Weise sind natürlich alle Terme außer dem ersten zu $0$ geworden; und es wird also
\[
	V\left( 1 + \frac{1}{2}z + \frac{1}{6} z^2 + \frac{1}{24}z^3 + \frac{1}{120}z^4 + \frac{1}{720}z^5 + \mathrm{etc}\right) = 1
\]
sein.
\paragraph{§7}
Weil also
\[
	e^z = 1 + z + \frac{z^2}{2} + \frac{z^3}{6} + \frac{z^4}{24} + \mathrm{etc}
\]
ist, wird
\[
	\frac{V(e^z -1)}{z} = 1
\]
sein und daher
\[
	V = \frac{z}{e^z - 1};
\]
damit dieser Ausdruck leichter wieder in eine Reihe verwandelt werden kann, wollen wir $z=2t$ setzen, dass
\[
	V = \frac{2t}{e^{2t} - 1}
\]
ist und daher
\[
	V + t = t\cdot\frac{e^{2t}+1}{e^{2t}-1}.
\]
Nun setze man
\[
	\frac{e^{2t}+1}{e^{2t}-1} = u
\]
und es wird
\[
	V = tu - t
\]
werden. Weil also
\[
	u = \frac{e^t + e^{-t}}{e^t - e^{-t}}
\]
ist, wird daher durch Entwicklung der Exponentialfunktionen
\[
	u = \frac{1 + \frac{1}{2}t^2 + \frac{1}{24}t^4 + \frac{1}{720}t^6 + \mathrm{etc.}}{t + \frac{1}{6}t^3 + \frac{1}{120}t^5 + \frac{1}{5040}t^7 + \mathrm{etc.}}
\]
sein, wo im Zähler allein die geraden Potenzen, im Nenner aber allein die ungeraden Potenzen auftauchen. Es ist aber klar, dass für ein sehr klein genommenes $t$ $u = \frac{1}{t}$ wird, die folgenden Terme aber nach den Potenzen $t$, $t^3$, $t^5$, etc fortschreiten werden.
\paragraph{§8}
Weil wir also
\[
	u = \frac{e^{2t} + 1}{e^{2t} - 1}
\]
haben, wird
\[
	e^{2t} = \frac{u+1}{u-1}
\]
sein und daher
\[
	2t = \log{\frac{u+1}{u-1}};
\]
Daher wird also durch Differenzieren
\[
	\partial t = -\frac{\partial u}{uu-1}
\]
sein, woher man folgert, dass
\[
	\frac{\partial u}{\partial t} + uu - 1 = 0
\]
sein wird. Weil wir aber wissen, dass der erste Term der Reihe, durch welche $u$ ausgedrückt wird, $\frac{1}{t}$ ist und die Potenzen der folgenden Potenzen um $2$ wachsen, setze man
\[
	u = \frac{1}{t} + 2At - 2Bt^3 + 2Ct^5 - 2Dt^7 + \mathrm{etc}
\]
und die Substitution geschehe auf folgende Weise:
\[
\begin{array}{rcccccccccccccc}
\frac{\partial u}{\partial t} &=& -\frac{1}{tt} &+& 2A &-& 6Bt^2 &+& 10ct^4 &-& 14Dt^6 &+& 18Et^8  &-& \mathrm{etc.}\\[1mm]
uu &=& +\frac{1}{tt} &+& 4A &-& 4B &+& 4C &-& 4D &+& 4E &-& \mathrm{etc.} \\[1mm]
& & & & &+& 4AA &-& 8AB &+& 8AC &-& 8AD &+& \mathrm{etc.} \\[1mm]
& & & & & & & & &+& 4BB &-& 8BC &+& \mathrm{etc.} \\[1mm]
-1 &=& -1 & & & & & & & & & & & &
\end{array}
\]
wo sich die ersten Terme von selbst aufheben, die übrigen aber die folgenden Bestimmungen liefern:
\begin{alignat*}{3}
6A  &= 1, &\qquad\text{also}\qquad & A = \tfrac{2}{3}\cdot\tfrac{1}{4} = \tfrac{1}{6} \\
10B &= 4AA, &\qquad\text{also}\qquad & B = \tfrac{2}{5}AA = \tfrac{1}{90} \\
14C &= 8AB, &\qquad\text{also}\qquad & C = \tfrac{2}{7}\cdot 2AB = \tfrac{1}{945} \\
18D &= 8AC+4BB, &\qquad\text{also}\qquad & D = \tfrac{2}{9}(2AC + BB) = \tfrac{1}{9540} \\
22E &= 8(AD+BC), &\qquad\text{also}\qquad & E = \tfrac{2}{11}(2AD + 2BC) = \tfrac{1}{93555} \\
\mathrm{etc.}
\end{alignat*}
\paragraph{§9}
Diese Buchstaben $A$, $B$, $C$, $D$, etc also sind genau dieselben, die ich einst benutzt habe um die Summen der Potenzen der Reziproken auszudrücken, weil ich ja gefunden habe, dass
\begin{align*}
	1 + \frac{1}{4} + \frac{1}{9} + \frac{1}{16} + \frac{1}{25} + \mathrm{etc} &= A\pi^2 \\
	1 + \frac{1}{4^2} + \frac{1}{9^2} + \frac{1}{16^2} + \frac{1}{25^2} + \mathrm{etc} &= B\pi^4 \\
	1 + \frac{1}{4^3} + \frac{1}{9^3} + \frac{1}{16^3} + \frac{1}{25^3} + \mathrm{etc} &= C\pi^6 \\
	\mathrm{etc}
\end{align*}
ist, welche Werte ich bis hin zur $34$.\,Potenz durch sehr arbeitsaufwendige Berechnungen ausgeführt habe.
\paragraph{§10}
Weil wir also
\[
	u = \frac{1}{t} + 2At - 2Bt^3 + 2Ct^5 - \mathrm{etc}
\]
genommen haben, wird wegen $V = tu-t$
\[
	V = 1 - t + 2At^2 - 2Bt^4 + 2Ct^6 - 2Dt^8 + \mathrm{etc}
\]
sein, wo nichts anderes übrig bleibt, außer dass $\frac{1}{2}z$ anstelle von $t$ geschrieben wird, woher
\[
	V = 1 - \frac{z}{2} + \frac{Azz}{2} - \frac{Bz^4}{8} + \frac{Cz^6}{32} - \frac{Dz^8}{128} + \mathrm{etc}
\]
hervorgeht. Weil wir daher
\[
	V = 1 + \alpha z + \beta z^2 + \gamma z^3 + \mathrm{etc}
\]
haben, werden wir nach der Anstellung des Vergleichs
\[
	\alpha = -\frac{1}{2},\quad \beta = \frac{1}{2}A,\quad \gamma = 0,\quad \delta = -\frac{1}{8}B,\quad \varepsilon = 0,\quad \zeta = \frac{1}{32}C,\quad \eta = 0,\quad \mathrm{etc}
\]
haben.
\paragraph{§11}
Nachdem die Werte dieser Buchstaben schon gefunden worden sind, wird die Summe der vorgelegten Reihe
\[
	S = X + X' + X'' + X''' + \mathrm{etc}
\]
auf folgende Weise ausgedrückt werden:
\[
	S  = -\int X\partial x + \frac{1}{2}X - \frac{1}{2}A\partial X + \frac{1}{8}B\partial^3 X - \frac{1}{32}C\partial^5 X + \frac{1}{128}D\partial^7 X - \frac{1}{512}E\partial^9 X + \mathrm{etc},
\]
wo das Integral $\int X\partial x$ so genommen werden muss, dass es für $x=\infty$ gesetzt verschwindet; daher ist klar, wenn die hinzuzufügende Konstante unendlich sein muss, dass dann auch die Summe selbst der Reihe unendlich sein muss.
\paragraph{§12}
Wir wollen das Beispiel betrachten, in welchem $X = \frac{1}{x^n}$ ist, sodass die Summe dieser Reihe zu sehen ist:
\[
	S = \frac{1}{x^n} + \frac{1}{(x+1)^n} + \frac{1}{(x+2)^n} + \frac{1}{(x+3)^n} + \mathrm{etc}.
\]
Hier wird also
\[
	\int X\partial x = -\frac{1}{(n-1)x^{n-1}}
\]
sein; damit diese Form für $x=\infty$ gesetzt verschwindet, ist es notwendig, dass der Exponent $n$ größer als die Einheit ist. Andernfalls würde nämlich, wenn $n=1$ oder $n<1$ wäre, die Summe der Reihe gewiss unendlich groß sein. Weiter wird in der Tat
\[
	\partial X = -\frac{n}{x^{n+1}}
\]
sein, daher
\[
	\partial^3 X = -\frac{n(n+1)(n+2)}{x^{n+3}},\quad \partial^5 X = -\frac{n(n+1)(n+2)(n+3)(n+4)}{x^{n+5}},\quad \mathrm{etc},
\]
nach Einsetzen welcher Werte die gesuchte Summe
\[
	S = \frac{1}{(n-1)x^{n-1}} + \frac{1}{2x^n} + \frac{A}{2}\cdot\frac{n}{x^{n+1}} - \frac{B}{8}\cdot\frac{n(n+1)(n+2)}{x^{n+3}} + \frac{C}{3^2}\cdot\frac{n\cdots (n+4)}{x^{n+5}} - \mathrm{etc}
\]
sein wird, welche Reihe umso mehr konvergieren wird, je größer die Zahl $x$ angenommen werden wird, außer dass die Buchstaben $A$, $B$, $C$, etc eine sehr konvergente Progression festsetzen.
\paragraph{§13}
Wenn also durch Beginnen von der Einheit aus diese Terme
\[
	1 + \frac{1}{2^n} + \frac{1}{3^n} + \frac{1}{4^n} + \cdots + \frac{1}{(x-1)^n}
\]
zusammengefasst werden und deren Summe $\triangle$ genannt wird, wird die Summe derselben ins Unendliche fortgesetzten Reihe $\triangle + S$ sein. Auf diese Weise habe ich einst die Summen solcher unendlichen Reihen für die einzelnen Werte $2$, $3$, $4$, $5$, etc des Exponenten $n$ auf viele Dezimalstellen berechnet, nachdem natürlich $x=10$ genommen worden ist, nach welcher Übereinkunft die Berechnung hinreichend angenehm ausgeführt werden könnte.
\section*{Fall 2 \\ Summation der unendlichen Reihe \\ $S = X - X' + X'' - X''' + X^4 - \mathrm{etc}$}
\paragraph{§14}
Wenn daher also der Index $x$ um die Einheit vermehrt wird, werden wir
\[
	S' = X' - X'' + X''' - X^4 + \mathrm{etc}
\]
haben. Man addiere diese Gleichung zu der vorhergehenden, und es wird diese endliche Gleichung hervorgehen
\[
	S + S' = X.
\]
Daher werden wir durch Differentialformeln
\[
	X = 2 S + \partial S + \frac{1}{2}\partial\partial S + \frac{1}{6}\partial^3 S + \frac{1}{2}\partial^4 S + \mathrm{etc}
\]
haben, woher nach Weglassen der Differentiale $S = \frac{1}{2}X$ sein wird, welcher also der erste Term der Reihe sein wird, die wir suchen. Wir wollen also
\[
	S = \frac{1}{2}X + \alpha \partial X + \beta\partial\partial X + \gamma \partial^3 X + \mathrm{etc}
\]
setzen und nach Ausführung der Substitution wird
\[
\begin{array}[t]{rccccccccccccc}
2S &= &X  &+& 2\alpha\partial X &+& 2\beta\partial\partial X &+& 2\gamma\partial^3 X &+& 2\delta\partial^4 X  &+& \mathrm{etc.} \\[1mm]
\partial S &= & & & \frac{1}{2}  &+& \alpha  &+& \beta  &+& \gamma &+& \mathrm{etc.} \\[1mm]
\frac{1}{2}\partial\partial &= & & & & & \frac{1}{4} &+& \frac{1}{2}\alpha &+& \frac{1}{2}\beta &+& \mathrm{etc.} \\[1mm]
\frac{1}{6}\partial^3 S &= & & & & & & & \frac{1}{12} &+& \frac{1}{6}\alpha &+& \mathrm{etc.} \\[1mm]
\frac{1}{24}\partial^4 S &= & & & & & & & & & \frac{1}{48} &+& \mathrm{etc.} \\[1mm]
\mathrm{etc}
\end{array}
\]
werden.
\paragraph{§15}
Nachdem also die ganzen Spalten $0$ gesetzt wurden, werden die folgenden Gleichheiten entstehen:
\[
	2\alpha + \frac{1}{2} = 0,\quad 2\beta + \alpha + \frac{1}{4} = 0,\quad 2\gamma + \beta + \frac{1}{2}\alpha + \frac{1}{12} = 0,
\]
\[ 
	2\delta + \gamma + \frac{1}{2}\beta + \frac{1}{6}\alpha 
	 +\frac{1}{48} = 0,\quad \mathrm{etc},
\]
woher zumindest die ersten Buchstaben diese Bestimmungen erhalten:
\[
	\alpha = -\frac{1}{4},\quad \beta = 0,\quad \gamma = \frac{1}{48},\quad \delta = 0,\quad \mathrm{etc.}
\]
\paragraph{§16}
Damit wir aber diese Werte leichter finden, wollen wir diese Reihe betrachten:
\[
	V = \frac{1}{2} + \alpha z + \beta z^2 + \gamma z^3 + \mathrm{etc},
\]
wessen Summe $V$ natürlich gesucht werden soll. Daher berechnen wir also die folgenden Reihen
\[
\begin{array}[t]{rccccccccccccccc}
2V &= &1  &+& 2\alpha z &+& 2\beta zz &+& 2\gamma z^3 &+& 2\delta z^4 &+& 2\varepsilon z^5 &+& \mathrm{etc.} \\[1mm]
Vz &= & & & \frac{1}{2}z &+& \alpha zz &+& \beta z^3  &+& \gamma z^4 &+& \delta z^5 &+&\mathrm{etc.} \\[1mm]
\frac{1}{2}Vzz &= & & & & & \frac{1}{4} &+& \frac{1}{2}\alpha &+& \frac{1}{2}\beta &+& \frac{1}{2}\gamma &+& \mathrm{etc.} \\[1mm]
\frac{1}{6}Vz^3 &= & & & & & & & \frac{1}{12} &+& \frac{1}{6}\alpha &+& \frac{1}{6}\beta &+& \mathrm{etc.} \\[1mm]
\frac{1}{24}Vz^4 &= & & & & & & & & & \frac{1}{48} &+& \frac{1}{24}\alpha &+& \mathrm{etc.} \\[1mm]
\mathrm{etc}
\end{array}
\]
Die Summe dieser Reihen wird also wegen der zuvor erwähnten Gleichheiten gleich $1$ werden, und so werden wir diese Gleichung haben:
\[
	V(2 + z + \frac{1}{2}z^2 + \frac{1}{6}z^3 + \frac{1}{24}z^4 + \mathrm{etc}) = 1.
\]
Weil daher
\[
	e^z = 1 + z + \frac{1}{2}z^2 + \frac{1}{6}z^3 + \mathrm{etc}
\]
ist, wird natürlich
\[
	V(1+e^z) = 1
\]
sein oder
\[
	V = \frac{1}{1+e^z},
\]
woher
\[
	2V - 1 = \frac{1 - e^z}{1 + e^z}.
\]
\paragraph{§17}
Man setze also wie zuvor
\[
	\frac{e^z - 1}{e^z + 1} = u,
\]
dass
\[
	2V = 1-u
\]
ist und es sei wiederum $2t=z$, sodass
\[
	u = \frac{e^t - e^{-t}}{e^t + e^{-t}}
\]
ist, und es wird nach Ausführung der Entwicklung
\[
	u = \frac{t + \frac{1}{6}t^3 + \frac{1}{120}t^5 + \frac{1}{5040}t^7 + \mathrm{etc}}{1 + \frac{1}{2}t^2 + \frac{1}{24}t^4 + \frac{1}{720}t^6 + \mathrm{etc}}
\]
sein. Daher ist klar, dass der erste Term dieser Reihe, der den Wert von u ausdrückt, $t$ sein wird, die folgenden aber durch ungerade Potenzen von $t$ fortschreiten.
\paragraph{§18}
Weil also 
\[
	u = \frac{e^{2t} - 1}{e^{2t}+1}
\]
ist, wird
\[
	e^{2t} = \frac{1+u}{1-u}
\]
sein und daher
\[
	2t = \log{\frac{1+u}{1-u}},
\]
woher durch Differenzieren
\[
\partial t = \frac{\partial u}{1-uu}
\]
wird, sodass
\[
	\frac{\partial u}{\partial t} + uu - 1 = 0
\]
ist, welches die Gleichung selbst ist, die für den ersten Fall gefunden worden ist. Und trotzdem geht deshalb für $u$ nicht dieselbe Reihe hervor. Weil ja hier nämlich der erste Term der Reihe gleich $t$ sein muss, ist eine Reihe dieser Art anzusetzen:
\[
	u = t - \mathfrak{A}t^3 + \mathfrak{B}t^5 - \mathfrak{C}t^7 + \mathfrak{D}t^9 - \mathfrak{E}t^{11} + \mathrm{etc},
\]
und es wird nach Ausführung der Substitution
\[
\begin{array}[t]{rccccccccccccccccc}
\frac{\partial u}{\partial t} &= &1 &-& 3\mathfrak{A}tt &+& 5\mathfrak{B}t^4 &-& 7\mathfrak{C}t^6 &+& 9\mathfrak{D}t^8 &-& 11\mathfrak{E}t^{10} &+& \mathrm{etc.} \\
uu &= & & & 1t &-& 2\mathfrak{A} &+& 2\mathfrak{B} &-& 2\mathfrak{C} &+& 2\mathfrak{D} &-& \mathrm{etc.} \\
& & & & & & &+& \mathfrak{A}^2 &-& 2\mathfrak{A}\mathfrak{B} &+& 2\mathfrak{A}\mathfrak{C} &-& \mathrm{etc.} \\
& & & & & & & & & & &-& \mathfrak{B}^2 &+& \mathrm{etc.} \\
-1 &= -&1  & & & & & & & & & & & &
\end{array}
\]
werden müssen und daher entstehen die folgenden Bestimmungen:
\[
\begin{array}[t]{rclcrcl}
\mathfrak{A} &=& 1 & \text{und daher} & \mathfrak{A} &=& \frac{1}{2} \\[1mm]
5\mathfrak{B} &=& 2\mathfrak{A} & \text{und daher} & \mathfrak{B} &=& \frac{2}{3}\mathfrak{A} = \frac{2}{15} \\[1mm]
7\mathfrak{C} &=& 2\mathfrak{B} + \mathfrak{A}^2 & \text{daher} & \mathfrak{C} &=& \frac{2}{7}\mathfrak{B} + \frac{1}{7}\mathfrak{A}^2 = \frac{17}{315} \\[1mm]
9\mathfrak{D} &=& 2\mathfrak{C} + 2\mathfrak{A}\mathfrak{B} & \text{also} & \mathfrak{D} &=& \frac{2}{9}\mathfrak{C} + \frac{2}{9}\mathfrak{A}\mathfrak{B} = \frac{62}{2835} \\[1mm]
\mathrm{etc.} & & & & \mathrm{etc.} & &
\end{array}
\]
\paragraph{§19}
Weil also
\[
	V = \frac{1}{2} - \frac{1}{2}u
\]
ist, wenn wir anstelle von $z$ wieder $\frac{z}{2}$ einsetzen, werden wir für $V$ diese Reihe finden:
\[
	V = \frac{1}{2} - \frac{1}{4}z + \frac{1}{16}\mathfrak{A}z^3 - \frac{1}{64}\mathfrak{B}z^5 + \frac{1}{256}\mathfrak{C}z^7 - \frac{1}{1024}\mathfrak{D}z^9 + \mathrm{etc}.
\]
Weil wir daher
\[
	V = \frac{1}{2} + \alpha z + \beta z^2 + \gamma z^3 + \delta z^4 + \mathrm{etc}
\]
gesetzt haben, berechnen wir daher die folgenden Werte der Buchstaben $\alpha$, $\beta$, $\gamma$, $\delta$, etc, die also
\begin{align*}
	\alpha &= -\frac{1}{4},\quad \beta = 0,\quad \gamma = \frac{1}{16}\mathfrak{A},\quad \delta = 0,\quad \varepsilon = -\frac{1}{64}\mathfrak{B},\quad \zeta = 0,\\
	\eta &= \frac{1}{256}\mathfrak{C},\quad \theta = 0
\end{align*}
sein werden; als logische Konsequenz wird die gesuchte Summe
\[
	S = \frac{1}{2}X - \frac{1}{4}\partial X + \frac{1}{16}\mathfrak{A}\partial^3 X - \frac{1}{64}\mathfrak{B}\partial^5 X + \frac{1}{256}\mathfrak{C}\partial^7 X - \mathrm{etc}
\]
sein.
\paragraph{§20}
Wir wollen nun diese Koeffizienten mit denen vergleichen, welche wir im vorhergehenden Fall für die einzelnen Differentiale erhalten haben, welche $\frac{A}{2}$, $\frac{B}{8}$, $\frac{C}{32}$, etc waren, und wir werden eine außerordentliche Relation zwischen jeden von beiden entdecken, wie sich aus diesem Schema sehen lässt:
\[
\begin{array}[t]{l|rcccl}
\partial X & \quad \frac{1}{4} : \frac{A}{2} &=& 3 &=& 2^2 - 1 \\[1mm]
\partial^3 X & \quad \frac{A}{16} : \frac{B}{8} &=& 15 &=& 2^4 - 1 \\[1mm]
\partial^5 X & \quad \frac{B}{64} : \frac{C}{32} &=& 63 &=& 2^6 - 1 \\[1mm]
\partial^7 X & \quad \frac{B}{256} : \frac{D}{128} &=& 255 &=& 2^8 - 1 \\[1mm]
\partial^9 X & \quad \frac{D}{1024} : \frac{E}{512} &=& 1023 &=& 2^{10} - 1 \\[1mm]
\mathrm{etc.} & \quad \mathrm{etc.} & & & & 
\end{array}
\]
\paragraph{§21}
Durch dieselben allbekannten Zahlen $A$, $B$, $C$, $D$, etc also wird auch in diesem Fall die gesuchte Summe auf die folgende Weise angenehm ausgedrückt werden:
\begin{align*}
	S =& \frac{1}{2}X - (2^2 - 1)\frac{1}{2}\cdot\partial X + (2^4 - 1)\frac{B}{8}\cdot\partial^3 X - (2^6 - 1)\frac{C}{32}\partial^5 X \\
	&+ (2^8 - 1)\frac{D}{128}\cdot\partial^7 X - (2^{10} - 1)\frac{E}{512}\cdot\partial^9 X + \mathrm{etc},
\end{align*}
welche Reihe sich, so weit es beliebt, fortsetzen lässt.
\end{document}